\documentclass[11pt,twoside, leqno]{article}

\usepackage{amssymb}
\usepackage{amsmath}
\usepackage{amsthm}
\usepackage[active]{srcltx}
\allowdisplaybreaks
\def\rr{{\mathbb R}}
\def\rn{{{\rr}^n}}
\def\rnz{{{\rr}^{n+1}_+}}

\def\zz{{\mathbb Z}}

\def\cc{{\mathbb C}}
\def\cn{{\mathbb N}}

\def\cs{{\mathcal S}}

\def\cm{{\mathcal M}}
\def\car{{\mathcal R}}

\def\ca{{\mathcal A}}
\def\ccc{{\mathcal C}}
\def\rz{{\rightarrow}}

\def\fz{\infty}
\def\az{\alpha}
\def\supp{{\mathop\mathrm{\,supp\,}}}
\def\dist{{\mathop\mathrm {\,dist\,}}}
\def\loc{{\mathop\mathrm{\,loc\,}}}
\def\lz{\lambda}
\def\dz{\delta}

\def\ez{\epsilon}

\def\bz{\beta}
\def\ro{\rho}

\def\gz{{\gamma}}
\def\oz{{\omega}}

\def\sz{\sigma}

\def\wz{\widetilde}
\def\nz{\nabla}

\def\hs{\hspace{0.3cm}}

\def\ls{\lesssim}

\def\bmom{{\mathrm{BMO}_{\ro,L}^{q,M}(\rn)}}
\def\bmo{{\mathrm{BMO}_{L}(\rn)}}
\def\bmor{{\mathrm{BMO}_{\ro,L}(\rn)}}
\def\bmop{{\mathrm{BMO}_{\ro,L^\ast}(\rn)}}
\def\cmo{{\mathrm{CMO}\,(\rn)}}

\def\lipc{{\mathrm{Lip}_{L^\ast}(\frac 1p -1,\rn)}}
\def\vmo{{\mathrm{VMO}_L(\rn)}}

\def\vmom{{\mathrm{VMO}_{\rho, L}^M(\rn)}}
\def\vmor{{\mathrm{VMO}_{\rho, L}(\rn)}}
\def\wvmo{{\wz{\mathrm{VMO}}_{\rho, L}^M(\rn)}}
\def\bw{{{B}_{\oz,L}(\rn)}}
\def\bp{{{B}_{\oz,L^\ast}(\rn)}}
\def\hw{{{H}_{\oz,L}(\rn)}}

\def\twl{{T_{\oz}({\rr}^{n+1}_+)}}

\def\twz{{ T^{\fz}_{\oz}({\rr}^{n+1}_+)}}
\def\tw0{{ T^{\fz}_{\oz,0}({\rr}^{n+1}_+)}}
\def\twv{{T^{\fz}_{\oz,\mathrm v}({\rr}^{n+1}_+)}}
\def\wtw{{\widetilde{T}_{\oz}({\rr}^{n+1}_+)}}

\def\com{\complement}
\def\r{\right}
\def\lf{\left}

\def\la{\langle}
\def\ra{\rangle}

\newtheorem{thm}{Theorem}[section]
\newtheorem{lem}{Lemma}[section]
\newtheorem{prop}{Proposition}[section]
\newtheorem{rem}{Remark}[section]
\newtheorem{cor}{Corollary}[section]
\newtheorem{defn}{Definition}[section]

\numberwithin{equation}{section}

\begin{document}

\arraycolsep=1pt
\title{{\vspace{-5cm}\small\hfill\bf Integral Equations Operator Theory, to appear}\\
\vspace{4cm}\Large\bf Generalized Vanishing Mean Oscillation Spaces
Associated with Divergence Form Elliptic Operators
\footnotetext{\hspace{-0.35cm}
2000 {\it Mathematics Subject Classification}. Primary 42B35;
Secondary 42B30, 46E30.
\endgraf{\it Key words and phrases.} divergence form
elliptic operator, Gaffney estimate, Orlicz function, Orlicz-Hardy space,
BMO, VMO, CMO, molecule, dual.
\endgraf
Dachun Yang is supported by the National Natural Science Foundation
(Grant No. 10871025) of China.
\endgraf $^\ast\,$Corresponding author.}}
\author{Renjin Jiang and Dachun Yang$\,^\ast$ }
\date{ }

\maketitle
\begin{center}
\begin{minipage}{13.5cm}\small
{\noindent{\bf Abstract.} Let $L$ be a divergence form elliptic
operator with complex bounded measurable coefficients,
$\omega$ a positive concave function on $(0,\infty)$
of strictly critical lower type $p_\omega\in (0, 1]$
and $\rho(t)={t^{-1}}/\omega^{-1}(t^{-1})$ for $t\in
(0,\infty).$ In this paper, the authors introduce
the generalized VMO spaces ${\mathop\mathrm{VMO}_
{\rho, L}({\mathbb R}^n)}$ associated with $L$,
and  characterize them via tent spaces. As applications,
the authors show that $(\mathrm{VMO}_{\rho,L}
({\mathbb R}^n))^\ast=B_{\omega,L^\ast}({\mathbb R}^n)$,
where $L^\ast$ denotes the adjoint operator of $L$
in $L^2({\mathbb R}^n)$ and $B_{\omega,L^\ast}({\mathbb R}^n)$ the Banach
completion of the Orlicz-Hardy space $H_{\omega,L^\ast}({\mathbb R}^n)$.
Notice that $\omega(t)=t^p$ for all $t\in (0,\infty)$ and
$p\in (0,1]$ is a typical example of positive concave functions
satisfying the assumptions. In particular, when $p=1$, then $\rho(t)\equiv 1$ and
$({\mathop\mathrm{VMO}_{1, L}({\mathbb R}^n)})^\ast=H_{L^\ast}^1({\mathbb R}^n)$,
where $H_{L^\ast}^1({\mathbb R}^n)$
was the Hardy space introduced by Hofmann and Mayboroda.}
\end{minipage}
\end{center}

\vspace{0.0cm}

\section{Introduction\label{s1}}
\hskip\parindent John and Nirenberg \cite{jn} introduced the space $\mathrm{BMO}(\rn)$,
which is defined to be
the space of all $f\in L^1_{\loc}(\rn)$ that satisfy
$$\|f\|_{\mathrm{BMO}(\rn)}\equiv \sup_{\mathrm{ball}\, B\subset \rn}\frac {1}{|B|}
\int_B \lf|f(x)-f_B\r|\,dx<\fz,$$
where and in what follows, $f_B\equiv \frac{1}{|B|}\int_B f(x)\,dx$.
The space $\mathrm{BMO}(\rn)$ is proved to be the dual of
the Hardy space $H^1(\rn)$ by Fefferman and Stein in \cite{fs}.

Sarason \cite{sa} introduced the space $\mathrm{VMO}(\rn)$,
which is defined to be
the space of all $f\in \mathrm{BMO}(\rn)$ that satisfy
$$\lim_{c\to 0}\sup_{B\subset \rn, \,r_B\le c}\frac {1}{|B|}
\int_B \lf|f(x)-f_B\r|\,dx=0.$$
In order to represent $H^1(\rn)$ as a dual space,
Coifman and Weiss \cite{cw} introduced the
space $\mathrm{CMO}(\rn)$ as the closure of continuous functions with
compact supports in the $\mathrm{BMO}(\rn)$ norm and showed that
$(\cmo)^\ast=H^1(\rn)$.
For more general cases, we refer to Janson \cite{ja}
and Bourdaud \cite{b}.

Let $L$ be a linear operator in $L^2(\rn)$ that generates an analytic semigroup
$\{e^{-tL}\}_{t\ge 0}$ with kernels satisfying an upper bound of
Poisson type. The Hardy space $H_L^1(\rn)$ and the BMO space
$\bmo$ associated with $L$ were defined and studied in \cite{adm,dy1}.
Duong and Yan \cite{dy2} further proved that $(H_L^1(\rn))^\ast=
\mathrm{BMO}_{L^\ast}(\rn)$, where and in what follows,
$L^\ast$ denotes the adjoint
operator of $L$ in $L^2(\rn).$
Moreover, recently, Deng et al in \cite{ddsty}
introduced the space $\mathrm{VMO}_L(\rn)$, the space of
vanishing mean oscillation associated with $L$, and proved that
$(\mathrm{VMO}_{L}(\rn))^\ast=H_{L^\ast}^1(\rn)$. Also,
Auscher and Russ \cite{ar} studied the Hardy space $H^1_L$ on
strongly Lipschitz domains associated with a divergence
form elliptic operator $L$ whose heat kernels have
the Gaussian upper bounds and regularity,
and Auscher, McIntosh and Russ \cite{amr}
treated the Hardy space $H^1$ associated with the Hodge
Laplacian on a Riemannian manifold with
doubling measure.

Let $A$ be an $n\times n$ matrix with entries
$\{a_{j,k}\}_{j,\,k=1}^n\subset L^\fz(\rn,\cc)$ satisfying the
uniform ellipticity condition, namely, there exist constants $0 <\lz_A\le
\Lambda_A<\fz$ such that for all $\xi,\,\zeta\in\cc^n$ and almost every $x\in\rn$,
\begin{equation}\label{1.1}
\lz_A|\xi|^2\le\car e \la A(x)\xi,\xi\ra \ \ \mbox{and} \ \ \ |\la
A(x)\xi,\zeta\ra| \le \Lambda_A |\xi||\zeta|.
\end{equation}
Then the second order divergence form operator is given by
\begin{equation}\label{1.2}
L f \equiv \mathop\mathrm{div}(A\nz f ),
\end{equation}
interpreted in the weak sense via a sesquilinear form. It is well
known that the kernel of the heat semigroup $\{e^{-tL}\}_{t>0}$
lacks pointwise estimates in general.

From now on, in this paper, we always let
$L$ be as in \eqref{1.2} and $L^\ast$ the adjoint operator
of $L$ in $L^2(\rn)$. Recently, Hofmann and Mayboroda
\cite{hm1, hm1c} studied the Hardy space
$H_L^1(\rn)$ and its dual space $\mathrm{BMO}_{L^\ast}(\rn)$.
Indeed, Hofmann and Mayboroda \cite{hm1} first defined
the Hardy space $H_L^1(\rn)$ via its molecular
decomposition, then established several maximal function characterizations of
$H_L^1(\rn)$, and in particular, showed that
$(H_L^1(\rn))^\ast=\mathrm{BMO}_{L^\ast}(\rn)$.
These results were generalized in \cite{jy2} to the Orlicz-Hardy
space $H_{\oz,L}(\rn)$ and its dual space
$\mathrm{BMO}_{\ro,L^\ast}(\rn)$, which contain
the Hardy spaces $H_L^p(\rn)$ for all $p\in (0,1]$,
the space $\mathrm{BMO}_{L^\ast}(\rn)$ and the Lipschitz spaces
$\lipc$ for all $p\in (0,1)$ as special cases.

Let $\omega$ be a positive concave function on $(0,\infty)$
of strictly critical lower type $p_\oz\in
(0, 1]$ and $\rho(t)={t^{-1}}/\omega^{-1}(t^{-1})$ for $t\in
(0,\infty).$ Recall that $\omega(t)=t^p$ for all $t\in (0,\fz)$ and
$p\in (0,1]$ is a typical example of such positive concave functions.
Motivated by \cite{cw,ja,ddsty,hm1,jy2}, in this paper,
we introduce the generalized VMO spaces $\vmor$
associated with $L$,
and  characterize them via tent spaces. Then, we prove that
$(\vmor)^\ast=\bp$, where $\bp$ denotes the Banach completion
of ${{H}_{\oz,L^\ast}(\rn)}$. When $\oz(t)=t$ for all
$t\in (0,\infty)$, we denote $\vmor$ simply by $\vmo$.
In this case, our result reads as $(\vmo)^\ast=H_{L^\ast}^1(\rn)$.
Finally, we show that the space $\cmo$
is a subspace of $\vmo$, and if $n\ge 3$, then there exists
an operator $L$ as in \eqref{1.2}, constructed in \cite{a1,f},
such that $\cmo\subsetneqq\vmo$. Moreover, when $n=1,\,2$,
the spaces $\cmo$ and $\vmo$ coincide with equivalent norms,
which is pointed to us by the referee.

Precisely, this paper is organized as follows. In Section \ref{s2}, we
recall some known definitions and notation on the divergence
form elliptic operator $L$ and Orlicz functions $\oz$ considered in this paper.

In Section \ref{s3}, we introduce the generalized VMO spaces $\vmor$
 associated with $L$, and tent spaces $\twv$ and give some basic properties of
these spaces. In particular, we
characterize the space $\vmor$ via $\twv$; see Theorem \ref{t3.1} below.

In Section \ref{s4}, we prove that
$(\mathrm{VMO}_{\rho,L}(\rn))^\ast=B_{\oz, L^\ast}(\rn)$,
where $B_{\oz, L^\ast}(\rn)$ denotes
the Banach completion of $H_{\oz, L^\ast}(\rn)$; see Theorem
\ref{t4.2} below. In particular, we have
$(\vmo)^\ast=H_{L^\ast}^1(\rn)$. Finally, in Proposition \ref{p4.2} below,
we show that the space $\cmo$ is a subspace of $\vmo$,
and if $n\ge 3$, then there exists an operator $L$ as in \eqref{1.2}
such that $\cmo\subsetneqq\vmo$; moreover, when $n=1,\,2$,
the spaces $\cmo$ and $\vmo$ coincide with equivalent norms.

Finally, we make some conventions. Throughout the paper, we denote
by $C$ a positive constant which is independent of the main
parameters, but it may vary from line to line. The symbol $X \ls Y$
means that there exists a positive constant $C$ such that $X \le
CY$; $B\equiv B(z_B,\,r_B)$ denotes an open ball with
center $z_B$ and radius $r_B$ and $CB(z_B,\,r_B)\equiv
B(z_B,\,Cr_B).$ Moreover, in what follows, for each ball
$B\subset\rn$ and $j\in\cn$,
we set $U_0(B)\equiv B$ and $U_j(B)\equiv2^jB\setminus 2^{j-1}B$.
Set $\cn\equiv\{1,2,\cdots\}$ and
$\zz_+\equiv\cn\cup\{0\}.$ For any subset $E$ of $\rn$, we denote by
$E^\com$ the set $\rn\setminus E.$

\section{Preliminaries\label{s2}}

\hskip\parindent  In this section, we recall some notions and notation on
divergence form elliptic operators, describe some basic assumptions on
Orlicz functions and also present some basic properties on them.

\subsection{Some notions on the divergence form elliptic operator $L$}

\hskip\parindent A family $\{S_t\}_{t>0}$ of operators  is said to satisfy the $L^2$
off-diagonal estimates, which are also called the Gaffney estimates
(see \cite{hm1}), if there exist positive constants $c,\,C$ and
$\bz$ such that for arbitrary closed sets $E,\,F\subset \rn$,
\begin{equation*}
   \|S_tf\|_{L^2(F)}\le Ce^{-(\frac{\dist(E,F)^2}{ct})^\bz}\|f\|_{L^2(E)}
\end{equation*}
for every $t>0$ and every $f\in L^2(\rn)$ supported in $E$. Here and
in what follows, for any $p\in (0,\fz]$ and $E\subset \rn$,
$\|f\|_{L^p(E)}\equiv \|f\chi_E\|_{L^p(\rn)}$; for any sets
$E,\,F\subset \rn$, $\dist(E,F)\equiv \inf\{|x-y|:\,x\in E,\,y\in
F\}.$

The following results were obtained in \cite{ahlmt,hm1,hm2}.

\begin{lem}[{\rm\cite{hm2}}]\label{l2.1} If two families,
$\{S_t\}_{t>0}$ and $\{T_t\}_{t>0}$, of operators
satisfy the Gaffney estimates, then so does
$\{S_tT_t\}_{t>0}$. Moreover, there exist positive constants
$c,\,C$ and $\bz$
 such that for arbitrary closed
sets $E,\,F \subset \rn$,
\begin{equation*}
\|S_sT_tf\|_{L^2(F)}\le Ce^{-(\frac{\dist(E,F)^2}
{c\max\{s,t\}})^\bz}\|f\|_{L^2(E)}
\end{equation*}
for every $s,\,t>0$ and every $f\in L^2(\rn)$ supported in $E$.
\end{lem}

\begin{lem}[{\rm\cite{ahlmt,hm2}}]\label{l2.2} The families
\begin{equation}\label{2.1}
 \{e^{-tL}\}_{t>0}, \ \ \{tLe^{-tL}\}_{t>0},
\end{equation}
as well as
\begin{equation}\label{2.2}
 \{(I+tL)^{-1}\}_{t>0},
\end{equation}
are bounded on $L^2(\rn)$ uniformly in $t$ and satisfy the Gaffney
estimates with positive constants $c$ and $C$,
depending on $n,\,\lz_A,\,\Lambda_A$ as in \eqref{1.1} only.
For the operators in \eqref{2.1}, $\bz=1$; while in \eqref{2.2}, $\bz=1/2$.
\end{lem}

Following
\cite{hm1}, set $
p_L\equiv \inf\{p\ge 1:\ \sup_{t>0}\|e^{-tL}\|_{L^p(\rn)\to L^p(\rn)}<\fz\}$
and
\begin{equation*}
 \wz p_L\equiv \sup\lf\{p\le \fz:\ \sup_{t>0}\|e^{-tL}\|_{L^p(\rn)
 \to L^p(\rn)}<\fz\r\}.
\end{equation*}
It was proved by Auscher \cite{a1} that
if $n=1,\,2$, then $p_L=1$ and $\wz p_L=\fz$,
and if $n\ge 3$, then $p_L<2n/(n+2)$ and $\wz p_L>2n/(n-2)$.
Moreover, thanks to a counterexample given by Frehse \cite{f},
this range is sharp.

\begin{lem}[{\rm\cite{hm1}}]\label{l2.3} Let $k\in\cn$ and
$p\in (p_L,\wz p_L)$. Then the operator given by setting,
for all $f\in L^p(\rn)$ and $x\in\rn$,
\begin{equation*}
\cs_L^kf(x)\equiv \bigg(\iint_{\Gamma(x)}|(t^2L)^ke^{-t^2L}f(y)|
\frac{\,dy\,dt}{t^{n+1}}\bigg)^{1/2},
\end{equation*}
is bounded on $L^p(\rn)$.
\end{lem}

\subsection{Orlicz functions \label{s2.2}}

\hskip\parindent Let $\omega$ be a positive function defined on
$\rr_+\equiv(0,\,\fz).$ The function $\omega$ is said to be of upper
type $p$ (resp. lower type $p$) for certain $p\in[0,\,\fz)$,  if there
exists a positive constant $C$ such that for all $t\geq 1$ (resp.
$t\in (0, 1]$) and all $s\in (0,\fz)$,
\begin{equation}\label{2.3}
\omega(st)\le Ct^p \omega(s).
\end{equation}

Obviously, if $\oz$ is of lower type $p$ for certain $p>0$, then
$\lim_{t\to0+}\oz(t)=0.$ So for the sake of convenience, if it is
necessary, we may assume that $\oz(0)=0.$ If $\oz$ is of both upper
type $p_1$ and lower type $p_0$, then $\oz$ is said to be of type
$(p_0,\,p_1).$ Let
\begin{equation*}
p_\oz^+\equiv\inf\{ p>0:\ \mathrm{there\ exists} \ C>0 \ \mathrm{such \ that }
\ \eqref{2.3} \ \mathrm{holds\ for\ all}\ t\in[1,\fz),\ s\in (0,\fz)\},
\end{equation*}
and
\begin{equation*}
p_\oz^-\equiv\sup\{ p>0:\ \mathrm{there\ exists} \ C>0 \ \mathrm{such \ that }
\ \eqref{2.3} \ \mathrm{holds\ for\ all}\
t\in(0,1],\ s\in (0,\fz)\}.
\end{equation*}
The function $\oz$ is said to be of strictly lower type $p$ if for all $t\in(0,1)$
and $s\in (0,\fz)$, $\omega(st)\le t^p \omega(s),$ and for a such function $\oz$,
we define
\begin{equation*}
p_\oz\equiv\sup\{ p>0: \omega(st)\le t^p \omega(s) \ \mathrm{holds\ for\
all}\ s\in (0,\fz)\ \mathrm{and}\ t\in(0,1)\}.
\end{equation*}
It is easy to see that $p_\oz\le p_\oz^{-}\le{p_\oz^+}$ for all $\oz.$
In what follows, $p_\oz$, $p_\oz^-$ and ${p_\oz^+}$ are called the
strictly critical lower type index, the critical lower
type index and the critical upper type index of $\oz$, respectively.

\begin{rem}\label{r2.1}\rm
We claim that if $p_\oz$ is defined as above,
then $\oz$ is also of strictly lower type $p_\oz$.
In other words, $p_\oz$ is attainable.
In fact, if this is not the case, then there exist
certain $s\in (0,\fz)$ and $t\in (0,1)$ such that
$\oz(st)>t^{p_\oz}\oz(s)$. Hence there exists $\ez\in(0,p_\oz)$
small enough such that $\oz(st)>t^{p_\oz-\ez }\oz(s)$, which is
contrary to the definition of $p_\oz$. Thus, $\oz$ is of strictly
lower type $p_\oz$.
\end{rem}

We now introduce the following assumption.

\begin{proof}[\bf Assumption (A)]\rm Let $\oz$ be a positive
function defined on $\rr_+$, which is of strictly lower type
and its strictly lower type index $p_\oz\in (0,1]$.
Also assume that $\oz$ is continuous, strictly
increasing and concave.
\end{proof}

Notice that if $\oz$ satisfies Assumption (A), then $\oz(0)=0$
and $\oz$ is obviously of upper type 1.
Since $\oz$ is concave,
it is subadditive. In fact, let $0<s<t$, then
$$\oz(s+t)\le \frac{s+t}{t}\oz(t)\le \oz(t)
+\frac{s}{t}\frac{t}{s}\oz(s)=\oz(s)+\oz(t).$$
For any concave function $\oz$ of strictly lower type $p$, if we set
$\wz\oz(t)\equiv\int_0^t\oz(s)/s\,ds$ for $t\in [0,\fz)$, then by
\cite[Proposition 3.1]{v}, $\wz\oz$ is equivalent to $\oz$, namely,
there exists a positive constant $C$ such that $C^{-1}\oz(t)\le
\wz\oz(t)\le C\oz(t)$ for all $t\in [0,\fz)$; moreover, $\wz\oz$
is strictly increasing, concave, subadditive and continuous function of
strictly lower type $p.$ Since all our results are invariant on equivalent
functions, we always assume that $\oz$ satisfies Assumption (A);
otherwise, we may replace $\oz$ by $\wz\oz.$

For example, if $\oz(t)=t^p$ with $p\in (0, 1]$, then $p_\oz=p_\oz^+=p$;
if $\oz(t)=t^{1/2}\ln(e^4+t)$, then $p_\oz=p_\oz^+=1/2$.

Let $\oz$ satisfy Assumption (A). A measurable function $f$ on
$\rn$ is said to be in the Lebesgue type space $L(\oz)$ if
$$\int_{\rn}\oz(|f(x)|)\,dx< \fz.$$
Moreover, for any $f\in L(\oz)$, define
$$\|f\|_{L(\oz)}\equiv\inf\lf\{\lz>0:\ \int_{\rn}\oz\lf(\frac{|f(x)|}
{\lz}\r)\,dx\le 1\r\}.$$

Since $\oz$ is strictly increasing, we define the function $\ro(t)$
on $\rr_+$ by setting, for all $t\in (0,\fz)$,
\begin{equation}\label{2.4}
\ro(t)\equiv\frac{t^{-1}}{\oz^{-1}(t^{-1})},
\end{equation}
where and in what follows, $\oz^{-1}$ denotes the inverse function of
$\oz.$ Then the types of $\oz$ and $\rho$ have the following
relation; see \cite{v} for its proof.

\begin{prop}\label{p2.1}
Let $0<p_0\le  p_1\le1$ and $\oz$ be an increasing function. Then $\oz$
is of type $(p_0,\, p_1)$ if and only if $\ro$ is of type
$(p_1^{-1}-1,\,p_0^{-1}-1).$
\end{prop}

Throughout the whole paper, we always assume that $\oz$ satisfies
Assumption (A) and $\ro$ is as in \eqref{2.4}.

\section{Spaces $\vmor$ associated with $L$\label{s3}}
\hskip\parindent In this section, we introduce the generalized vanishing
mean oscillation spaces associated with $L$.
We begin with some notions and notation.

 Let $q\in (p_L, \wz p_L)$, $M\in\cn$ and $\ez\in(0,\fz)$.
 A function $\az\in L^q(\rn)$ is called an $(\oz,q,M,\ez)_L$-molecule
adapted to $B$ if there exists a ball $B$ such that

{\rm (i)} $\|\az\|_{L^q(U_j(B))}\le 2^{-j\ez}|2^jB|^{1/q-1}\ro(|2^jB|)^{-1}$,
$j\in {\zz}_+$;

{\rm (ii)} For every $k=1,\cdots,M$ and $j\in {\zz}_+$, there holds
$$\|(r_B^{-2}L^{-1})^{k}\az\|_{L^q(U_j(B))}\le
2^{-j\ez}|2^jB|^{1/q-1}\ro(|2^jB|)^{-1}.$$

Let $\ez$ and $M$ be as above. We also introduce the space
\begin{equation}\label{3.1}
\cm_\oz^{M,\ez}(L)\equiv \{\mu\in L^2(\rn):\,
\|\mu\|_{\cm_\oz^{M,\ez}(L)}<\fz\},\end{equation}
where
$$\|\mu\|_{\cm_\oz^{M,\ez}(L)}\equiv \sup_{j\ge 0}
\lf\{2^{j\ez}|B(0,2^j)|^{1/2}
\ro(|B(0,2^j)|)\sum_{k=0}^M\|L^{-k}\mu\|_{L^2(U_j(B(0,1)))}\r\}.$$

Notice that if $\phi\in \cm_\oz^{M,\ez}(L)$ with norm 1, then $\phi$ is an
$(\oz,2,M,\ez)$-molecule adapted to $B(0,1)$. Conversely, if $\az$
is an $(\oz,2,M,\ez)$-molecule adapted to certain ball, then $\az\in
\cm_\oz^{M,\ez}(L)$.

Let $A_t$ denote either $(I+t^2L)^{-1}$ or $e^{-t^2L}$ and $f\in
(\cm_\oz^{M,\ez}(L))^\ast$, the dual of $\cm_\oz^{M,\ez}(L)$. We
claim that $(I-A_t^\ast)^Mf\in L^2_{\loc}(\rn)$ in the sense of distributions.
 In fact, for any ball $B$, if $\psi\in L^2(B)$, then it follows from the Gaffney
estimates via Lemmas \ref{l2.1} and \ref{l2.2} that $(I-A_t)^M\psi\in
\cm_\oz^{M,\ez}(L)$ for all $\ez>0$ and any fixed $t\in(0,\fz)$. Thus,
\begin{eqnarray*}
  |\la (I-A_t^\ast)^Mf,\psi\ra|\equiv|\la f,(I-A_t)^M\psi\ra|\le
  C(t,r_B,\dist(B,0))\|f\|_{(\cm_\oz^{M,\ez}(L))^\ast}\|\psi\|_{L^2(B)},
\end{eqnarray*}
which implies that $(I-A_t^\ast)^Mf\in L^2_{\loc}(\rn)$  in the sense of distributions.

Finally, for any $M\in \cn$, define
\begin{equation}\label{3.2}
\cm_{\oz,L^\ast}^M(\rn)\equiv \bigcap_{\ez>n(1/p_\oz-1/p_\oz^+)}
(\cm_\oz^{M,\ez}(L))^\ast.\end{equation}

\begin{defn}\label{d3.1}
Let $q\in(p_L,\wz p_L)$ and  $M> \frac n2 (\frac1{p_\oz}-\frac 12)$.
 An element  $f\in\cm_{\oz,L}^M(\rn)$ is said to be in $\bmom$ if
\begin{equation*}
\|f\|_{\bmom}\equiv\sup_{B\subset\rn}\frac{1}{\ro(|B|)}\lf[\frac{1}{|B|}\int_B
|(I-e^{-r_B^2L})^Mf(x)|^q \,dx\r]^{1/q}< \fz,
\end{equation*}
where the supremum is taken over all balls $B$ of $\rn.$
\end{defn}

\begin{rem}\rm
The spaces $\bmom$ were introduced in \cite{jy2};
moreover, if $\oz(t)=t$ for all $t\in (0,\fz)$,
$\bmom$ is the space $\bmo$ introduced by Hofmann
and Mayboroda \cite{hm1}. Since the spaces $\bmom$
coincide for all $q\in(p_L,\wz p_L)$ and
$M> \frac n2 (\frac 1{p_\oz}-\frac 12)$ (see \cite[Theorem 4.1]{jy2}),
in what follows, we denote $\bmom$
simply by $\bmor$.
\end{rem}
Let us introduce a new space $\vmom$ as a subspace of $\bmor$.
\begin{defn}\label{d3.2}
Let $M> \frac n2 (\frac1{p_\oz}-\frac 12)$.
 An element $f\in\bmor$ is said to be in $\vmom$ if it satisfies
 the following limiting conditions $\gz_1(f)=\gz_2(f)=\gz_3(f)=0$, where
$$\gz_1(f)\equiv\lim_{c\to 0}\sup_{B:\,r_B\le c}\frac{1}{\ro(|B|)}\bigg(
\frac{1}{|B|}\int_B|(I-e^{-r_B^2L})^Mf(x)|^2\,dx\bigg)^{1/2},$$
$$\gz_2(f)\equiv\lim_{c\to \fz}\sup_{B:\,r_B\ge c}\frac{1}{\ro(|B|)}\bigg(
\frac{1}{|B|}\int_B|(I-e^{-r_B^2L})^Mf(x)|^2\,dx\bigg)^{1/2},$$
and
$$\quad\gz_3(f)\equiv\lim_{c\to \fz}\sup_{B\subset B(0,c)^\com}
\frac{1}{\ro(|B|)}\bigg(
\frac{1}{|B|}\int_B|(I-e^{-r_B^2L})^Mf(x)|^2\,dx\bigg)^{1/2}.$$
 For any $f\in\vmom$, define $\|f\|_{\vmom}\equiv\|f\|_{\bmor}$.
\end{defn}

We next present some properties of the space $\vmom$. To this end, we
first recall some notions of tent spaces; see
\cite {cms}.

Let $\Gamma(x)\equiv\{(y,t)\in\rnz:\,|x-y|<t\}$ denote the standard
cone (of aperture 1) with vertex $x\in\rn.$ For any closed set $F$
of $\rn$, denote by $\car{F}$ the union of all cones with vertices
in $F$, i.\,e., $\car{F}\equiv\cup_{x\in F}\Gamma(x)$; for any
open set $O$ in $\rn$, denote the tent over $O$ by $\widehat{O}$,
which is defined by $\widehat{O}\equiv[\car(O^\com)]^\com.$

For any measurable function $g$ on ${\rr}^{n+1}_+$ and $x\in\rn$, define
\begin{equation*}
\ca(g)(x)\equiv
\lf(\iint_{\Gamma(x)}|g(y,t)|^2\frac{\,dy\,dt}{t^{n+1}}\r)^{1/2}
\end{equation*}
and
\begin{equation*}
\ccc_{\ro}(g)(x)\equiv\sup_{\mathrm{ball}\, B\ni x}\frac{1}{\ro(|B|)}
\lf(\frac{1}{|B|}\int_{\widehat{B}}|g(y,t)|^2\frac{\,dy\,dt}{t}\r)^{1/2}.
\end{equation*}

For $p\in(0,\fz)$, the tent space $T^p_2(\rnz)$ is defined to be the space
of all measurable functions $g$ on $\rnz$ such that $\|g\|_{T^p_2(\rnz)}\equiv
\|\ca(g)\|_{L^p(\rn)}<\fz.$
Let $\oz$ satisfy Assumption (A). The
tent space $T_\oz(\rnz)$ associated to the function $\oz$
is defined to be the space of measurable functions $g$ on $\rr^{n+1}_+$ such that
$\ca(g)\in L(\oz)$ with the norm defined by
$$\|g\|_{T_\oz(\rnz)}\equiv \|\ca(g)\|_{L(\oz)}=\inf\lf\{\lz>0:\ \int_{\rn} \oz
\lf(\frac{\ca(g)(x)}{\lz}\r) \,dx\le 1\r\};$$
the space $T^{\fz}_{\oz}({\rr}^{n+1}_+)$  is defined to be the
space of all measurable functions $g$ on $\rr^{n+1}_+$ satisfying
$\|g\|_{T^{\fz}_{\oz}({\rr}^{n+1}_+)}\equiv\|\ccc_{\ro}(g)\|_{L^\fz(\rn)}<\fz$.

In what follows,  let $\twv$ be the
space of all $f\in\twz$ satisfying $\eta_1(f)=\eta_2(f)=\eta_3(f)=0$,
where
\begin{equation*}
\eta_1(f)\equiv \lim_{c\to0}\sup_{B:\,r_B\le c}\frac{1}{\ro(|B|)}
\lf(\frac{1}{|B|}\iint_{\widehat{B}}|f(y,t)|^2\frac{\,dy\,dt}{t}\r)^{1/2},
\end{equation*}
$$\eta_2(f)\equiv \lim_{c\to\fz}\sup_{B:\,r_B\ge c}\frac{1}{\ro(|B|)}
\lf(\frac{1}{|B|}\iint_{\widehat{B}}|f(y,t)|^2\frac{\,dy\,dt}{t}\r)^{1/2},$$
and
$$\quad \eta_3(f)\equiv \lim_{c\to\fz}\sup_{B\subset B(0,c)^\com}
\frac{1}{\ro(|B|)}
\lf(\frac{1}{|B|}\iint_{\widehat{B}}|f(y,t)|^2\frac{\,dy\,dt}{t}\r)^{1/2}.$$
It is easy to see that $\twv$ is a closed linear subspace of $\twz$.

Further, denote by $T^\fz_{\oz,1}(\rnz)$ the space of all $f\in\twz$
with $\eta_1(f)=0$, and $T^2_{2,c}(\rnz)$ the space of all $f\in
T^2_{2}(\rnz)$ with compact supports. Obviously, we have
$T^2_{2,c}(\rnz)\subset \twv\subset T^\fz_{\oz,1}(\rnz)$.
Finally, denote by $\tw0$ the closure of
$T_{2,c}^2(\rnz)$ in the space $T^\fz_{\oz,1}(\rnz)$.

By \cite[Proposition 3.1]{jy1}, we have the following result; see also \cite{ddsty}.

\begin{lem}\label{l3.1}
Let $\twv$ and $\tw0$ be defined as above. Then $\twv=\tw0$.
\end{lem}

Recall that a measure $d\mu$ on ${\rr}^{n+1}_+$ is called a
$\ro$-Carleson measure if
\begin{equation*}
\sup_{B\subset
\rn}\frac{1}{\ro(|B|)}\lf\{\frac{1}{|B|}\iint_{\widehat{B}}|d\mu|\r\}^{1/2}<\fz,
\end{equation*}
where the supremum is taken over all balls $B$ of $\rn$.

We now characterize the space $\vmom$ via tent spaces.

\begin{thm}\label{t3.1}
Let $M,M_1\in\cn$ and
$M_1\ge M>\frac n2(\frac 1{p_\oz}-\frac12).$
Then the following conditions are equivalent:

(a) $f\in \vmom$;

(b) $f\in \cm_{\oz,L}^{M_1}(\rn)$ and $(t^2L)^{M_1}e^{-t^2L}f\in \twv$.

Moreover, $\|(t^2L)^{M_1}e^{-t^2L}f\|_{\twz}$ is equivalent to
$\|f\|_{\bmor}.$
\end{thm}

To prove Theorem \ref{t3.1}, we need two auxiliary results.

Let $M\in\zz_+$. In what follows, let $C_M$ be the positive constant such that
\begin{equation}\label{3.3}
C_{M}\int^\fz_0 t^{2(M+1)}e^{-2t^2}\frac{\,dt}{t}=1.
\end{equation}
The following lemma was established in \cite[Proposition 4.6]{jy2}.

\begin{lem}\label{l3.2}
Let $\ez\in(0,\fz)$ and $M>\frac n2(\frac 1{p_\oz}-\frac12).$
If $f\in\bmor$, then for any $(\oz,2,M,\ez)_{L^\ast}$-molecule $\az$,
there holds
\begin{equation*}
\int_{\rn}f(x)\az(x)\,dx=C_{M}\iint_{\rnz}
(t^2L)^{M} e^{-t^2L}f(x) t^2L^\ast e^{-t^2L^\ast}\az(x)\,\frac{dx\,dt}{t}.
\end{equation*}
\end{lem}

\begin{defn}\label{d3.3}
Let $M>\frac n2(\frac 1{p_\oz}-\frac12).$
Define the space $\wvmo$ to be the space of all elements
$f\in\bmor$ that satisfy the limiting
conditions $\wz\gz_1(f)=\wz\gz_2(f)=\wz\gz_3(f)=0$, where
$$\wz\gz_1(f)\equiv\lim_{c\to 0}\sup_{B:\,r_B\le c}\frac{1}{\ro(|B|)}\bigg(
\frac{1}{|B|}\int_B|(I-[I+r_B^2L]^{-1})^Mf(x)|^2\,dx\bigg)^{1/2},$$
$$\wz\gz_2(f)\equiv\lim_{c\to \fz}\sup_{B:\,r_B\ge c}\frac{1}{\ro(|B|)}\bigg(
\frac{1}{|B|}\int_B|(I-[I+r_B^2L]^{-1})^Mf(x)|^2\,dx\bigg)^{1/2},$$
and
$$\quad\wz\gz_3(f)\equiv\lim_{c\to \fz}\sup_{B\subset B(0,c)^\com}
\frac{1}{\ro(|B|)}\bigg(
\frac{1}{|B|}\int_B|(I-[I+r_B^2L]^{-1})^Mf(x)|^2\,dx\bigg)^{1/2}.$$
\end{defn}

\begin{prop}\label{p3.1}
Let $p\in(0,1]$ and $M>\frac n2(\frac 1p-\frac12).$
Then $f\in\vmom$ if and only if $f\in\wvmo$
\end{prop}
\begin{proof}\rm
Recall that it was proved in \cite[Lemma 4.1]{jy2} that
\begin{eqnarray}\label{3.4}
\|f\|_{\bmor}&&\sim\sup_{B\subset\rn}\frac{1}{\ro(|B|)}\lf[\frac{1}{|B|}\int_B
|(I-(I+r_B^2L)^{-1})^Mf(x)|^2 \,dx\r]^{1/2}.
\end{eqnarray}

Now suppose that $f\in\wvmo$. To see $f\in\vmom$,
it suffices to show that
\begin{eqnarray}\label{3.5}
  &&\frac{1}{\ro(|B|)|B|^{1/2}}\bigg(\int_B|(I-e^{-r_B^2L})^Mf(x)|^2
  \,dx\bigg)^{1/2}\ls\sum_{k=0}^\fz  2^{-k}\dz_k(f,B),
\end{eqnarray}
where
\begin{eqnarray}\label{3.6}
\dz_k(f,B)\equiv&& \sup_{\{B'\subset
2^{k+1}B:\ r_{B'}\in[2^{-1}r_B,r_B]\}}\frac{1}{\ro(|B|)|B|^{1/2}}\\
&&\hspace{2cm}\times\bigg\{\int_{B'}|(I-[I+r_B^2L]^{-1})^Mf(x)|^2
\,dx\bigg\}^{1/2}.\nonumber
\end{eqnarray}

Indeed, since $f\in\wvmo$, by Definition \ref{d3.3} and \eqref{3.4},
we have that $\dz_k(f,B)\ls
\|f\|_{\bmor}$ and for each $k\in\cn$,
\begin{eqnarray*}
 \lim_{c\rz 0}\sup_{B:\,r_B\le c}\dz_k(f,B)&&=\lim_{c\rz
 \fz}\sup_{B:\,r_B\ge c}\dz_k(f,B)=\lim_{c\rz \fz}\sup_{B\subset
 B(0,\,c)^\com}\dz_k(f,B)=0.
\end{eqnarray*}
Then by the dominated convergence
theorem for series, we have
\begin{eqnarray*}
\gz_1(f)&&=\lim_{c\rz 0}\sup_{B:\,r_B\le c}
\frac{1}{\ro(|B|)|B|^{1/2}}\bigg(\int_{B}|(I-e^{-r_B^2L})^Mf(x)|^2
\,dx\bigg)^{1/2}\\
&&\ls \sum_{k=1}^{\fz}2^{-k}\lim_{c\rz
0}\sup_{B:\,r_B\le c}\dz_k(f,B)=0.\nonumber
\end{eqnarray*}
Similarly, we have that $\gz_2(f)=\gz_3(f)=0$, and thus, $f\in\vmom$.

Let us now prove \eqref{3.5}. Write
\begin{equation}\label{3.7}
f=(I-[I+r_B^2L]^{-1})^Mf+\{I-(I-[I+r_B^2L]^{-1})^M\}f\equiv f_1+f_2,
\end{equation}
where by Lemma \ref{l2.2}, we have
\begin{eqnarray}\label{3.8}
  &&\|(I-e^{-r_B^2L})^Mf_1\|_{L^2(B)}\\
  &&\hs\le\sum_{k=0}^\fz
  \|(I-e^{-r_B^2L})^M(f_1\chi_{U_k(B)})\|_{L^2(B)}\ls\sum_{k=0}^\fz
  e^{-c2^{2k}}\|f_1\chi_{U_k(B)}\|_{L^2(\rn)}\nonumber\\
  &&\hs\ls\ro(|B|)|B|^{1/2}\sum_{k=0}^\fz
  e^{-c2^{2k}}2^{kn}\dz_k(f,B) \ls\ro(|B|)|B|^{1/2}\sum_{k=0}^\fz
  2^{-k}\dz_k(f,B),\nonumber
\end{eqnarray}
where $c$ is a positive constant and the third inequality follows from the fact that
there exists a collection, $\{B_{k,1},B_{k,2},\cdots, B_{k,N_k}\}$, of balls
such that each ball $B_{k,i}$ is
of radius $r_{B}$, $B(x_B, 2^{k+1}r_B)\subset\cup^{N_k}_{i=1} B_{k,i}$
and $N_k \ls 2^{nk}.$

To estimate the remaining term, by the formula that
\begin{equation}\label{3.9}
I-(I-[I+r_B^2L]^{-1})^M=
\sum_{j=1}^M \frac{M!}{j!(M-j)!} (r_B^2L)^{-j}(I-[I+r_B^2L]^{-1})^M
\end{equation}
(which relies on the fact that $(I-(I+r^2L)^{-1})(r^2L)^{-1}=(I+r^2L)^{-1}$
for all $r\in (0,\fz)$),
and the Minkowski inequality, we obtain
\begin{eqnarray}\label{3.10}
\quad\quad&&\|(I-e^{-r_B^2L})^Mf_2\|_{L^2(B)}\\
 &&\hs\ls \sum_{j=1}^M
 \bigg(\int_B\lf|(I-e^{-r_B^2L})^{M-j}\bigg[-\int_0^{r_B}
 \frac{s}{r_B^2}e^{-s^2L}\,ds\bigg]^jf_1(x)\r|^2\,dx\bigg)^{1/2}\nonumber\\
 &&\hs\ls\sum_{j=1}^M\sum_{k=0}^{M-j}
 \int_0^{r_B}\cdots\int_0^{r_B}\frac{s_1}{r_B^2}\cdots\frac{s_j}{r_B^2}
 \|e^{-(kr_B^2+s_1^2+\cdots+s_j^2)L}
 f_1\|_{L^2(B)}\,ds_1\cdots\,ds_j\nonumber\\
&&\hs\ls \sum_{j=1}^M\sum_{k=0}^{M-j}\int_0^{r_B}\cdots
\int_0^{r_B}\frac{s_1}{r_B^2}\cdots\frac{s_j}{r_B^2}
 \sum_{i=0}^\fz e^{-\frac{c(2^ir_B)^2}{kr_B^2+s_1^2+\cdots+s_j^2}}
 \|f_1\chi_{U_i(B)}\|_{L^2(\rn)}\,ds_1\cdots\,ds_j\nonumber\\
 &&\hs\ls \ro(|B|)|B|^{1/2}\sum_{i=0}^\fz e^{-\frac{c2^{2i}}{M}}2^{in}
 \dz_i(f,B)\ls \ro(|B|)|B|^{1/2}\sum_{i=0}^\fz 2^{-i}
 \dz_i(f,B),\nonumber
 \end{eqnarray}
where $c$ is a positive constant and in the penultimate inequality,
we used the fact that $\int_0^{r_B}\cdots
\int_0^{r_B}\frac{s_1}{r_B^2}\cdots\frac{s_j}{r_B^2}\,ds_1\cdots\,ds_j
\sim 1$. Combining the estimates \eqref{3.8}
and \eqref{3.10}, we obtain \eqref{3.5},
which further implies that  $\wvmo\subset\vmom$.

By borrowing some ideas from the proof of \cite[Lemma 8.1]{hm1},
then similarly to the proof above, we have that
$\vmom\subset \wvmo$. We omit the details here, which completes
the proof of Proposition \ref{p3.1}.
\end{proof}

\begin{proof}[Proof of Theorem \ref{t3.1}]
We first show that (a) implies (b). Let
$f\in \vmom$. By \cite[Theorem 6.1]{jy2},
we have that $(t^2L)^{M_1}e^{-t^2L}f\in \twz$. To see
 $(t^2L)^{M_1}e^{-t^2L}f\in \twv$, we claim that it suffices to show that
for all balls $B$,
\begin{eqnarray}\label{3.11}
&&\frac{1}{\ro(|B|)|B|^{1/2}}\bigg(\iint_{\widehat B}|(t^2L)^{M_1}e^{-t^2L}f(x)|^2
\frac{\,dx\,dt}{t}\bigg)^{1/2}\ls\sum_{k=0}^\fz  2^{-k}\dz_k(f,B),
\end{eqnarray}
where $\dz_k(f,B)$ is as in \eqref{3.6}.
In fact, since $f\in\vmom=\wvmo$, we have that for each $k\in\cn$,
$\dz_k(f,B)\ls \|f\|_{\bmor}$ and
\begin{eqnarray*}
 \lim_{c\rz 0}\sup_{B:\,r_B\le c}\dz_k(f,B)&&=\lim_{c\rz
 \fz}\sup_{B:\,r_B\ge c}\dz_k(f,B)=\lim_{c\rz \fz}\sup_{B\subset
 B(0,\,c)^\com}\dz_k(f,B)=0.
\end{eqnarray*}
Then by the dominated convergence
theorem for series, we have
\begin{eqnarray*}
\eta_1(f)&&=\lim_{c\rz 0}\sup_{B:\,r_B\le c}
\frac{1}{\ro(|B|)|B|^{1/2}}\bigg(\iint_{\widehat B}|(t^2L)^{M_1}e^{-t^2L}f(x)|^2
\frac{\,dx\,dt}{t}\bigg)^{1/2}\\
&&\ls \sum_{k=1}^{\fz}2^{-k}\lim_{c\rz
0}\sup_{B:\,r_B\le c}\dz_k(f,B)=0.\nonumber
\end{eqnarray*}
Similarly, we have $\eta_2(f)=\eta_3(f)=0$, and thus,
$(t^2L)^{M_1}e^{-t^2L}f\in\twv$.

Let us prove \eqref{3.11}. Write $f\equiv f_1+f_2$ as in \eqref{3.7}.
Then by Lemmas \ref{l2.2} and \ref{l2.3}, similarly to the estimate \eqref{3.8},
we obtain
\begin{eqnarray}\label{3.12}
  &&\bigg(\iint_{\widehat B}|(t^2L)^{M_1}e^{-t^2L}f_1(x)|^2
  \frac{\,dx\,dt}{t}\bigg)^{1/2}\\
  &&\hs\le \sum_{k=0}^\fz
  \bigg(\iint_{\widehat B}|(t^2L)^{M_1}e^{-t^2L}(f_1\chi_{U_k(B)})(x)|^2
  \frac{\,dx\,dt}{t}\bigg)^{1/2}\nonumber\\
  &&\hs\ls \|f_1\|_{L^2(4B)}+  \sum_{k=3}^\fz\bigg(\int_0^{r_B}
\exp\lf\{-\frac{(2^kr_B)^2}{ct^2}\r\}
  \frac{\,dt}{t}\bigg)^{1/2}\|f_1\chi_{U_k(B)}\|_{L^2(\rn)}\nonumber\\
  &&\hs\ls  \|f_1\|_{L^2(4B)}+  \sum_{k=3}^\fz\bigg\{\int_0^{r_B}
\lf[\frac{t^2}{(2^kr_B)^2}\r]^{n+1}
  \frac{\,dt}{t}\bigg\}^{1/2}\|f_1\chi_{U_k(B)}\|_{L^2(\rn)}\nonumber\\
&&\hs\ls \ro(|B|)|B|^{1/2}\sum_{k=0}^{\fz}2^{-k}\dz_k(f,B),\nonumber
\end{eqnarray}
where $c$ is a positive constant.
Applying \eqref{3.9}, Lemma \ref{l2.2} and $M_1\ge M$ to $f_2$ yields that
\begin{eqnarray}\label{3.13}
  &&\bigg(\iint_{\widehat B}|(t^2L)^{M_1}e^{-t^2L}f_2(x)|^2
  \frac{\,dx\,dt}{t}\bigg)^{1/2}\\
  &&\hs\ls \sum_{j=1}^{M}\bigg(\iint_{\widehat B}|(t^2L)^{M_1}
e^{-t^2L}(r_B^2L)^{-j}f_1(x)|^2\frac{\,dx\,dt}{t}\bigg)^{1/2}\nonumber\\
  &&\hs \ls\sum_{j=1}^M\sum_{k=0}^\fz
  \bigg(\int_{\widehat B}\lf(\frac{t^2}{r_B^2}\r)^{2j}|(t^2L)^{M_1-j}
e^{-t^2L}(f_1\chi_{U_k(B)})(x)|^2\frac{\,dx\,dt}{t}\bigg)^{1/2}\nonumber\\
&&\hs\ls\sum_{j=1}^M\lf\{\sum_{k=0}^2
  \bigg[\int_0^{r_B}\lf(\frac{t^2}{r_B^2}\r)^{2j}
  \frac{\,dt}{t}\bigg]^{1/2}\r.\nonumber\\
&&\hs\hs\lf.+\sum_{j=1}^M\sum_{k=3}^\fz
  \bigg[\int_0^{r_B}\lf(\frac{t^2}{r_B^2}\r)^{2j}\exp\lf\{-\frac{(2^kr_B)^2}{ct^2}\r\}
  \frac{\,dt}{t}\bigg]^{1/2}\r\}\|f_1\chi_{U_k(B)}\|_{L^2(\rn)}\nonumber\\
&&\hs\ls  \|f_1\|_{L^2(4B)}+  \sum_{k=3}^\fz\bigg[\int_0^{r_B}
\lf(\frac{t^2}{(2^kr_B)^2}\r)^{n+1}
  \frac{\,dt}{t}\bigg]^{1/2}\|f_1\chi_{U_k(B)}\|_{L^2(\rn)}\nonumber\\
&&\hs\ls \ro(|B|)|B|^{1/2}\sum_{k=0}^{\fz}2^{-k}\dz_k(f,B).\nonumber
\end{eqnarray}

The estimates \eqref{3.12} and \eqref{3.13} imply \eqref{3.11},
and hence, completes the proof that (a) implies (b).

Conversely, let $f\in \cm_{\oz,L}^{M_1}(\rn)$ and
$(t^2L)^{M_1}e^{-t^2L}f\in \twv$.
By \cite[Theorem 6.1]{jy2}, we have that $f\in\bmor$. For any ball $B$,
write
\begin{eqnarray*}
  \bigg(\int_{B}|(I-e^{-r_B^2L})^Mf(x)|^2
\,dx\bigg)^{1/2}&&=\sup_{\|g\|_{L^2(B)}\le 1}
\lf|\int_{B}(I-e^{-r_B^2L})^Mf(x)g(x)\,dx\r|\\
&&=\sup_{\|g\|_{L^2(B)}\le 1}
\lf|\int_{\rn}f(x)(I-e^{-r_B^2L^\ast})^M g(x)\,dx\r|.
\end{eqnarray*}
Notice that for any $g\in L^2(B)$,
$(I-e^{-r_B^2L^\ast})^M g$ is
a multiple of an $(\oz,2,M,\ez)_{L^\ast}$-molecule;
see \cite[p.\,43]{hm1} and \cite{jy2}.
Then by Lemma \ref{l3.2} and the H\"older inequality, we obtain
\begin{eqnarray*}
&& \bigg(\int_{B}|(I-e^{-r_B^2L})^Mf(x)|^2
\,dx\bigg)^{1/2}\\
&&\hs =\sup_{\|g\|_{L^2(B)}\le 1}
\lf|C_{M_1}\iint_{\rnz}
(t^2L)^{M_1} e^{-t^2L}f(x) t^2L^\ast e^{-t^2L^\ast}
(I-e^{-r_B^2L^\ast})^M g(x)\frac{\,dx\,dt}{t}\r|\\
&&\hs\ls \sum_{k=0}^\fz\lf\{\iint_{V_k(B)}|(t^2L)^{M_1}
e^{-t^2L}f(x)|^2\frac{\,dx\,dt}{t}\r\}^{1/2}\\
&&\hs\hs\times
\sup_{\|g\|_{L^2(B)}\le 1}\lf\{\iint_{V_k(B)}|(t^2L^\ast e^{-t^2L^\ast}
(I-e^{-r_B^2L^\ast})^M g(x)|^2\frac{\,dx\,dt}{t}\r\}^{1/2}
\equiv \sum_{k=0}^\fz \sz_k(f,B)\mathrm{I}_k,
\end{eqnarray*}
where $V_0(B)\equiv\widehat B$ and $V_k(B)\equiv(\widehat {2^kB})
\setminus (\widehat {2^{k-1}B})$ for $k\in\cn$.
In what follows, for $k\ge 2$, let $V_{k,1}(B)\equiv (\widehat {2^kB})
\setminus ({2^{k-2}B}\times (0,\fz))$ and
$V_{k,2}(B)\equiv V_k(B)\setminus V_{k,1}(B)$.

For $k=0,1,2$, by Lemmas \ref{2.2} and \ref{l2.3}, we obtain
\begin{eqnarray*}
\mathrm{I}_k&&=\sup_{\|g\|_{L^2(B)}\le 1}\lf\{\iint_{V_k(B)}|(t^2L^\ast e^{-t^2L^\ast}
(I-e^{-r_B^2L^\ast})^M g(x)|^2\frac{\,dx\,dt}{t}\r\}^{1/2}\\
&&\ls  \sup_{\|g\|_{L^2(B)}\le 1}
\|(I-e^{-r_B^2L^\ast})^M g\|_{L^2(\rn)}\ls  1.
\end{eqnarray*}
Now for $k\ge 3$, write
\begin{eqnarray*}
\mathrm{I}_k&&\ls \sup_{\|g\|_{L^2(B)}\le 1}\lf\{\iint_{V_{k,1}(B)}
|(t^2L^\ast e^{-t^2L^\ast}
(I-e^{-r_B^2L^\ast})^M g(x)|^2\frac{\,dx\,dt}{t}\r\}^{1/2}\\
&&\hs +\sup_{\|g\|_{L^2(B)}\le 1}\lf\{\iint_{V_{k,2}(B)}
\cdots\r\}^{1/2}\equiv \mathrm{I}_{k,1}+\mathrm{I}_{k,2}.
\end{eqnarray*}
Since for any $(y,t)\in V_{k,2}(B)$, $t\ge 2^{k-2}r_B$,
by the Minkowski inequality and Lemma \ref{l2.2}, we obtain
\begin{eqnarray*}
\mathrm{I}_{k,2}&&\ls \sup_{\|g\|_{L^2(B)}\le 1}\lf\{\iint_{V_{k,2}(B)}
\lf|t^2L^\ast e^{-t^2L^\ast}\bigg[-\int_0^{r_B^2}
 L^\ast e^{-sL^\ast}\,ds\bigg]^Mg(x)\r|^2\frac{\,dx\,dt}{t}\r\}^{1/2}\\
&&\ls \sup_{\|g\|_{L^2(B)}\le 1}\int_0^{r_B^2}\cdots\int_0^{r_B^2}
\lf\{\iint_{V_{k,2}(B)}
|t^2(L^\ast)^{M+1} \r.\\
&&\hs\hs \lf.\times e^{-(t^2+s_1+\cdots+s_M)L^\ast}g(x)|^2
\frac{\,dx\,dt}{t}\r\}^{1/2}\,ds_1\cdots\,ds_M\\
&&\ls\sup_{\|g\|_{L^2(B)}\le 1}\int_0^{r_B^2}\cdots\int_0^{r_B^2}
\lf\{\int_{2^{k-2}r_B}^{2^kr_B} \frac{t^4\|g\|^2_{L^2(B)}}
{(t^2+s_1+\cdots+s_M)^{2(M+1)}}
\frac{\,dt}{t}\r\}^{1/2}\,ds_1\cdots\,ds_M\\
&&\ls 2^{-2kM}.
\end{eqnarray*}
Similarly,  we have that $\mathrm{I}_{k,1}\ls 2^{-2kM}$. Combining the
above estimates and the fact that $\ro$ is of upper type
$1/p_\oz-1$, we finally obtain
\begin{eqnarray*}
\frac{1}{\ro(|B|)|B|^{1/2}}\bigg(\int_{B}|(I-e^{-r_B^2L})^Mf(x)|^2
\,dx\bigg)^{1/2}&&\ls \sum_{k=0}^\fz 2^{-2kM}
\frac{1}{\ro(|B|)|B|^{1/2}}\sz_k(f,B)\\
&&\ls \sum_{k=0}^\fz 2^{-k[2M-n(\frac 1{p_\oz} -\frac 12)]}
\frac{\sz_k(f,B)}{\ro(|2^kB|)|2^kB|^{1/2}}.
\end{eqnarray*}

Since $(t^2L)^{M_1}e^{-t^2L}f\in\twv\subset \twz$, by  $M>\frac n2
(\frac 1{p_\oz} -\frac 12)$ and the dominated convergence theorem for
series, we have
\begin{eqnarray*}
\gz_1(f)&&=\lim_{c\rz 0}\sup_{B:\,r_B\le c}
\frac{1}{\ro(|B|)|B|^{1/2}}\bigg(\int_{B}|(I-e^{-r_B^2L})^Mf(x)|^2
\,dx\bigg)^{1/2}\\
&&\ls \sum_{k=0}^\fz 2^{-k[2M-n(\frac 1{p_\oz} -\frac 12)]}
\lim_{c\rz 0}\sup_{B:\,r_B\le c}
\frac{\sz_k(f,B)}{\ro(|2^kB|)|2^kB|^{1/2}}=0.
\end{eqnarray*}
Similarly, $\gz_2(f)=\gz_3(f)=0$,
which implies that $f\in \vmom$,
and hence, completes the proof of Theorem \ref{t3.1}.
\end{proof}

\begin{rem}\rm
  (i) It follows from Theorem \ref{t3.1} that for all $M\in\cn$
  and $M> \frac n2(\frac 1{p_\oz}-\frac 12)$, the spaces $\vmom$
  coincide with equivalent norms. Thus, in what follows,
  we denote $\vmom$ simply by $\vmor$; in particular, if
  $\oz(t)\equiv t$ for all $t\in (0,\fz)$, then $\ro(t)\equiv 1$ and
 we denote $\vmom$  simply by ${\mathop\mathrm{VMO}_L}(\rn)$.

 (ii) When $A$ has real entries, or when the dimension $n=1$ or $n=2$ in
the case of complex entries, the heat kernels always satisfy the
Gaussian pointwise estimates (see \cite{at}), in these cases,
the space ${\mathop\mathrm{VMO}_L}(\rn)$ here coincides with the one
introduced by Deng et al in \cite{ddsty}.
\end{rem}

\section{Dual spaces of $\vmor$ \label{s4}}
\hskip\parindent In this section, we identify the dual spaces of
$\vmor$. We begin with some notions and known
facts on tent spaces.

Recall that a function $a$ on $\rr^{n+1}_+$ is called a
$T_{\oz}({\rr}^{n+1}_+)$-atom if there exists a ball $B\subset
\rn$ such that $\supp a\subset \widehat{B}$ and
$(\iint_{\widehat{B}}|a(x,t)|^2\frac{\,dx\,dt}{t})^{1/2}\le
\frac{1}{\ro(|B|)|B|^{1/2}}.$

\begin{defn}\label{d4.1} Let $p\in (0,1)$.
The space $\wtw$ is defined to be the space of
all $f=\sum_{j\in\cn}\lz_j a_j$ in $(\twz)^\ast$, where $\{a_j\}_{j\in\cn}$ are
$\twl$-atoms and $\{\lz_j\}_{j\in\cn}\in \ell^1.$ If $f\in\wtw$, then define
$\|f\|_{\wtw}\equiv\inf\{\sum_{j\in\cn}|\lz_j|\},$
where the infimum is taken over all the possible decompositions of $f$ as above.
\end{defn}

By \cite[Lemma 3.1]{hm1}, $\wtw$ is a Banach space. Moreover,
from Definition \ref{d4.1}, it is easy to deduce that $\twl$ is dense in
$\wtw$; in other words, $\wtw$ is a Banach completion of $\twl$.
The following lemma is just \cite[Lemma 4.1]{jy1}.

 \begin{lem}\label{l4.1}
  $\twl\subset \wtw$, $\twl$ is dense in $\wtw$
 and there exists a positive constant
 $C$ such that for all $f\in\twl$, $ \|f\|_{\wtw}\le C\|f\|_{\twl}.$
\end{lem}

The following theorem was established in \cite{wws}; see also \cite[Theorem 4.2]{jy1}.
\begin{thm}\label{t4.1}
$(\twv)^\ast=\wtw$.
\end{thm}

Now, let us recall some notions on the Hardy spaces
associated with $L$.
For all $f\in L^2(\rn)$ and $x\in\rn$, define
\begin{equation*}
  \cs_Lf(x)\equiv \bigg(\iint_{\Gamma(x)}|t^2Le^{-t^2L}f(y)|^2
  \frac{\,dy\,dt}{t^{n+1}}\bigg)^{1/2}.
  \end{equation*}
The space $\hw$ is defined to be the completion of the
set $\{ f\in L^2(\rn):\, \cs_Lf\in L(\oz)\}$ with respect to the
quasi-norm  $\|f\|_{\hw}\equiv \|\cs_Lf\|_{L(\oz)}$.

The Orlicz-Hardy space $H_{\oz,L}(\rn)$ was introduced
and studied in \cite{jy2}. If $\oz(t)\equiv t$ for all
$t\in (0,\fz)$, then the space
$\hw$ coincides with the Hardy space $H_L^1(\rn)$,
which was introduced and studied by
Hofmann and Mayboroda \cite{hm1} (see also \cite{hm1c}).

\begin{defn}\label{d4.2} Let $M\in\cn$, $M>\frac n2
(\frac1{p_\oz}-\frac 12)$ and
$\ez\in (n(\frac1{p_\oz}-\frac 1{p_\oz^+}),\fz)$.
An element $f\in (\bmop)^\ast$ is said to be in
the space $H_{\oz, L}^{M,\ez}(\rn)$
 if there exist $\{\lz_j\}_{j=1}^\fz\subset \cc$ and $(\oz,2,M,\ez)_L$-molecules
  $\{\az_j\}_{j=1}^\fz$ such that
$f=\sum_{j=1}^\fz\lz_j\az_j$ in $(\bmop)^\ast$ and
$$\Lambda(\{\lz_j\az_j\}_j)\equiv\inf\lf\{\lz>0:\, \sum_{j=1}^\fz|B_j|\oz
\lf(\frac{|\lz_j|}{\lz|B_j|\ro(|B_j|)}\r)\le1\r\}<\fz,$$
where for each $j$, $\az_j$ is adapted to the ball $B_j$.

If $f\in H_{\oz, L}^{M,\ez}(\rn)$, then its norm is defined by
$\|f\|_{H_{\oz, L}^{M,\ez}(\rn)}\equiv \inf\Lambda(\{\lz_j\az_j\}_j),$
where the infimum is taken over all the possible decompositions of $f$ as above.
\end{defn}
It was proved in \cite[Theorem 5.1]{jy2} that for all
$M>\frac n2 (\frac1{p_\oz}-\frac 12)$ and
$\ez\in (n(\frac1{p_\oz}-\frac 1{p_\oz^+}),\fz)$,
the spaces $\hw$ and $H_{\oz,L}^{M,\ez}(\rn)$
coincide with equivalent norms.

Let us introduce the Banach completion of the space $H_L^p(\rn)$.

\begin{defn}\label{d4.3} Let $\ez\in(n(\frac1{p_\oz}-\frac 1{p_\oz^+}),\fz)$
 and $M>\frac n2 (\frac 1{p_\oz}-\frac 12)$.
The space $B_{\oz,L}^{M,\ez}(\rn)$ is defined to be
the space of all $f=\sum_{j\in\cn}\lz_j\az_j$ in
$(\bmop)^\ast$, where  $\{\lz_j\}_{j\in\cn}\in
\ell^1$ and $\{\az_j\}_{j\in\cn}$ are $(\oz,2,M,\ez)_L$-molecules.
If $f\in B_{\oz,L}^{M,\ez}(\rn)$, define
$\|f\|_{B_{\oz,L}^{M,\ez}(\rn)}\equiv\inf\{\sum_{j\in\cn}|\lz_j|\},$
where the infimum is taken over all the possible decompositions of $f$
as above.
\end{defn}

By \cite[Lemma 3.1]{hm1} again, we see that $B_{\oz,L}^{M,\ez}(\rn)$ is
a Banach space. Moreover, from Definition \ref{d4.2}, it is easy to deduce that
$\hw$ is dense in $B_{\oz,L}^{M,\ez}(\rn)$.
More precisely, we have the following lemma.

\begin{lem}\label{l4.2} Let $\ez\in(n(\frac1{p_\oz}-\frac 1{p_\oz^+}),\fz)$
and $M>\frac n2 (\frac 1{p_\oz}-\frac 12)$. Then

(i) $\hw\subset B_{\oz,L}^{M,\ez}(\rn)$ and the inclusion is
continuous.

(ii) For any $\ez_1\in(n(\frac1{p_\oz}-\frac 1{p_\oz^+}),\fz)$
and $M_1>\frac n2 (\frac 1{p_\oz}-\frac 12)$,
the spaces $B_{\oz,L}^{M,\ez}(\rn)$ and $B_{\oz,L}^{M_1,\ez_1}(\rn)$
coincide with equivalent norms.
\end{lem}

\begin{proof}\rm From Definition \ref{d4.3} and the molecular characterization
of $\hw$, it is easy
to deduce (i). Let us prove (ii). By symmetry, it suffices
to show that $B_{\oz,L}^{M,\ez}(\rn)\subset B_{\oz,L}^{M_1,\ez_1}(\rn)$.

Let $f\in B_{\oz,L}^{M,\ez}(\rn)$. By Definition \ref{d4.3},
there exist $(\oz,2,M,\ez)_L$-molecules
$\{\az_j\}_{j\in\cn}$ and $\{\lz_j\}_{j\in\cn}\subset\cc$ such that
$f=\sum_{j\in \cn}\lz_j\az_j$ in $(\bmop)^\ast$ and $\sum_{j\in\cn}|\lz_j|
\ls\|f\|_{B_{\oz,L}^{M,\ez}(\rn)}$.

By (i), for each $j\in\cn$, we have that $\az_j\in \hw\subset
B_{\oz,L}^{M_1,\ez_1}(\rn)$ and
$\|\az_j\|_{B_{\oz,L}^{M_1,\ez_1}(\rn)}\ls \|\az_j\|_{\hw}\ls 1$.
Since $B_{\oz,L}^{M_1,\ez_1}(\rn)$ is a Banach space, we obtain that
$f\in B_{\oz,L}^{M_1,\ez_1}(\rn)$ and
$\|f\|_{B_{\oz,L}^{M_1,\ez_1}(\rn)}\le
\sum_{j\in\cn}|\lz_j|\|\az_j\|_{B_{L,M_1}^{p,\ez_1}(\rn)}
\ls \|f\|_{B_{\oz,L}^{M,\ez}(\rn)}$. Thus,
$B_{\oz,L}^{M,\ez}(\rn)\subset B_{\oz,L}^{M_1,\ez_1}(\rn)$,
which completes the proof of Lemma \ref{l4.2}.
\end{proof}

Since the spaces $B_{\oz,L}^{M,\ez}(\rn)$ coincide for all
$\ez\in(n(\frac1{p_\oz}-\frac 1{p_\oz^+}),\fz)$
and $M>\frac n2 (\frac 1{p_\oz}-\frac 12)$,
in what follows, we denote $B_{\oz,L}^{M,\ez}(\rn)$ simply by
$B_{\oz,L}(\rn)$.

\begin{lem}\label{l4.3} $(B_{\oz,L})^\ast=\bmop.$
\end{lem}

\begin{proof}
Since $(\hw)^\ast=\bmop$ and $\hw\subset \bw$,
by duality, we have that $(\bw)^\ast\subset\bmop.$

Conversely, let $\ez\in(n(\frac1{p_\oz}-\frac 1{p_\oz^+}),\fz)$,
$M>\frac n2 (\frac 1{p_\oz}-\frac 12)$
and $f\in\bmop$. For any $g\in \bw$, by Definition \ref{d4.3},
there exist $\{\lz_j\}_{j\in\cn}\subset \cc$ and
$(\oz,M,2,\ez)_L$-molecules $\{\az_j\}_{j\in\cn}$
such that $g=\sum_{j\in\cn}\lz_j\az_j$ in
$(\bmop)^\ast$ and $\sum_{j\in\cn}|\lz_j|\le 2\|g\|_{\bw}$.
Thus,
 \begin{eqnarray*}|\la f,g \ra|&&\le \sum_{j\in\cn}|\lz_j||\la f,\az_j\ra|
\ls \sum_{j\in\cn}|\lz_j|\|f\|_{\bmop}\|\az_j\|_{\hw}\\
&&\ls \|f\|_{\bmop}\|g\|_{\bw},
\end{eqnarray*}
which implies that $f\in (\bw)^\ast$,
and hence, completes the proof of Lemma \ref{l4.3}.
\end{proof}

Let $M\in\cn$. For all $F\in L^2(\rnz)$ with compact support, define
\begin{equation}\label{4.1}
  \pi_{L,M}F\equiv C_M\int_0^\fz (t^2L)^{M}e^{-t^2L}F(\cdot,t)\frac{\,dt}{t},
\end{equation}
where $C_M$ is the positive constant same as in
\eqref{3.3}.

\begin{prop}\label{p4.1}
Let $M\in\cn$. Then the operator $\pi_{L,M}$, initially
defined on $T_{2,c}^2(\rnz)$, extends to a bounded linear operator

{\rm (i)} from $T_2^{p}(\rnz)$ to $L^p(\rn),$ if $p\in (p_L,\wz p_L)$;

{\rm (ii)} from $\twl$ to $\hw$,  if $M>\frac n2 (\frac1{p_\oz}-\frac 12)$;

{\rm (iii)} from $\wtw$ to $\bw$, if $M>\frac n2 (\frac1{p_\oz}-\frac 12)$;

{\rm (iv)} from $\twv$ to $\vmor$.
\end{prop}

\begin{proof} (i) and (ii) were established in \cite{jy2}.
Let us prove (iii).

By \cite[Lemma 4.7]{jy1}, we have that $T_{2,c}^2(\rnz)$ is dense in $\wtw$.
Let $f\in T_{2,c}^2(\rnz)$. By (ii) and Lemma \ref{l4.2}, we obtain that
$\pi_{L,M}f\in \hw\subset \bw$. Moreover,
by Definition \ref{d4.1}, there exist $\twl$-atoms $\{a_j\}_{j\in\cn}$
and $\{\lz_j\}_{j\in\cn}\subset\cc$ such that
$f=\sum_{j\in\cn} \lz_ja_j$ in $(\twz)^\ast$ and
$\sum_j|\lz_j|\ls\|f\|_{\wtw}.$
In addition, for any $g\in\bmop$, by \cite[Theorem 6.1]{jy2},
we have that $(t^2L^\ast)^{M}e^{-t^2L^\ast}g\in \twz$. Thus, by the fact
that $(\twl)^\ast=\twz$ (see \cite[Proposition 4.1]{jyz}), we obtain
\begin{eqnarray*}
\la\pi_{L,M}(f),g\ra&&=C_M\int_{\rnz}f(x,t)
(t^2L^\ast)^{M}e^{-t^2L^\ast}g(x)\,\frac{dx\,dt}{t}\\
&&=\sum_{j\in\cn}\lz_jC_M\int_{\rnz}a_j(x,t)(t^2L^\ast)^{M}
e^{-t^2L^\ast}g(x)\,\frac{dx\,dt}{t}=\sum_{j\in\cn}\lz_j\la\pi_{L,M}(a_j),g\ra,
\end{eqnarray*}
which implies that $\pi_{L,M}(f)=\sum_{j\in\cn}\lz_j\pi_{L,M}(a_j)$
 in $(\bmop)^\ast$. By (ii), we further obtain
\begin{eqnarray*}\|\pi_{L,M}(f)\|_{\bw}&&\le \sum_{j\in\cn}|\lz_j|
 \|\pi_{L,M}(a_j)\|_{\bw}\\
 &&\ls\sum_j|\lz_j|\|\pi_{L,M}(a_j)\|_{\hw} \ls \|f\|_{\wtw}.
 \end{eqnarray*}
Since $T_{2,c}^2(\rnz)$ is dense in $\wtw$, we see that $\pi_{L,M}$ extends
to a bounded linear operator from $\wtw$ to $\bw$, which
completes the proof of (iii).

Let us now prove (iv). By the definition of $\tw0$, we have that
$T_{2,c}^2(\rnz)$ is dense in $\tw0$. Then by Lemma
\ref{l3.1}, we obtain that $T_{2,c}^2(\rnz)$ is dense
in $\twv$. Thus, to prove (iv), it suffices to show that $\pi_{L,M}$
maps $T_{2,c}^2(\rnz)$ continuously into $\vmor$.

Let $f\in T_{2,c}^2(\rnz)$. By (i), we see that
$\pi_{L,M}f\in L^2(\rn)$. Notice that \eqref{3.1} and \eqref{3.2}
with $L$ and $L^\ast$ exchanged implies that $L^2(\rn)\subset \cm_{\oz,L}^{M_1}(\rn)$,
where $M_1\in\cn$ and $M_1>\frac n2 (\frac 1{p_\oz} -\frac 12)$.
Thus, $\pi_{L,M}f\in\cm_{\oz,L}^{M_1}(\rn)$. To show $\pi_{L,M}f\in \vmor$,
by Theorem \ref{t3.1}, we only need to verify that
$(t^2L)^{M_1}e^{-t^2L}\pi_{L,M}f\in \twv$.

For any ball $B\equiv B(x_B,r_B)$, let
$V_0(B)\equiv\widehat B$ and $V_k(B)\equiv(\widehat {2^kB})
\setminus (\widehat {2^{k-1}B})$ for any $k\in\cn$.
For all $k\in\zz_+$, let $f_k\equiv f\chi_{V_k(B)}$.
Thus, for $k=0,1,2$, by Lemma \ref{l2.3} and (i), we obtain
\begin{eqnarray*}
  &&\lf(\iint_{\widehat B} |(t^2L)^{M_1}e^{-t^2L}
\pi_{L,M}f_k(x)|^2\frac{\,dx\,dt}{t}\r)^{1/2}\ls \|\pi_{L,M}f_k\|_{L^2(\rn)}
\ls \|f_k\|_{T_2^2(\rnz)} .
\end{eqnarray*}
For $k\ge 3$, let $V_{k,1}(B)\equiv (\widehat {2^kB})
\setminus ({2^{k-2}B}\times (0,\fz))$ and
$V_{k,2}(B)\equiv V_k(B)\setminus V_{k,1}(B)$. We further
write $f_k=f_k\chi_{V_{k,1}(B)}+f_k\chi_{V_{k,2}(B)}\equiv f_{k,1}+f_{k,2}$.
By the Minkowski inequality, Lemma \ref{l2.2}
and the H\"older inequality, we obtain
\begin{eqnarray*}
  &&\lf(\iint_{\widehat B} |(t^2L)^{M_1}e^{-t^2L}
\pi_{L,M}f_{k,2}(x)|^2\frac{\,dx\,dt}{t}\r)^{1/2}\\&&
\hs= C_M\lf(\iint_{\widehat B} \lf|\int_{2^{k-2}r_B}^{2^kr_B}(t^2L)^{M_1}e^{-t^2L}
(s^2L)^{M}e^{-s^2L}(f_{k,2}(\cdot,s))(x)\frac{\,ds}{s}\r|^2
\frac{\,dx\,dt}{t}\r)^{1/2}\\
&&\hs\ls \int_{2^{k-2}r_B}^{2^kr_B}
\lf(\iint_{\widehat B} \lf|t^{2M_1}s^{2M}
L^{M+M_1}e^{-(s^2+t^2)L}
(f_{k,2}(\cdot,s))(x)\r|^2
\frac{\,dx\,dt}{t}\r)^{1/2}\frac{\,ds}{s}\\
&&\hs\ls \int_{2^{k-2}r_B}^{2^kr_B}
\lf(\int_0^{r_B} \lf|\frac {t^{2M_1}s^{2M}}{(s^2+t^2)^{M+M_1}}\r|^2
\|f_{k,2}(\cdot,s)\|^2_{L^2(\rn)}
\frac{\,dt}{t}\r)^{1/2}\frac{\,ds}{s}\\
&&\hs\ls 2^{-2kM_1}\int_{2^{k-2}r_B}^{2^kr_B}
\|f_{k,2}(\cdot,s)\|_{L^2(\rn)}\frac{\,ds}{s}\ls
2^{-2kM_1}\|f_{k,2}\|_{T_2^2(\rnz)}.
\end{eqnarray*}
Similarly, we have
\begin{eqnarray*}
  &&\lf(\iint_{\widehat B} |(t^2L)^{M_1}e^{-t^2L}
\pi_{L,M}f_{k,1}(x)|^2\frac{\,dx\,dt}{t}\r)^{1/2}\ls
2^{-2kM_1}\|f_{k,1}\|_{T_2^2(\rnz)}.
\end{eqnarray*}

 Combining the above estimates, we finally obtain that
\begin{eqnarray*}
  &&\frac{1}{\ro(|B|)|B|^{1/2}}\lf(\iint_{\widehat B} |(t^2L)^{M_1}e^{-t^2L}
\pi_{L,M}f(x)|^2\frac{\,dx\,dt}{t}\r)^{1/2}\\
&&\hs\ls\sum_{k=0}^2 \frac{1}{\ro(|B|)|B|^{1/2}}
\lf(\iint_{\widehat B} |(t^2L)^{M_1}e^{-t^2L}
\pi_{L,M}f_k(x)|^2\frac{\,dx\,dt}{t}\r)^{1/2}\\
&&\hs\hs+\sum_{k=3}^\fz \sum_{i=1}^2\frac{1}{\ro(|B|)|B|^{1/2}}\lf(\iint_{\widehat B}
|(t^2L)^{M_1}e^{-t^2L}
\pi_{L,M}f_{k,i}(x)|^2\frac{\,dx\,dt}{t}\r)^{1/2}\\
&&\hs\ls\sum_{k=0}^2\frac{1}{\ro(|B|)|B|^{1/2}}\|f_k\|_{T_2^2(\rnz)}
+\sum_{k=3}^\fz \sum_{i=1,2}\frac{2^{-2kM_1}}{\ro(|B|)|B|^{1/2}}
\|f_{k,i}\|_{T_2^2(\rnz)}\\
&&\hs\ls\sum_{k=0}^\fz 2^{-k[2M_1-n(\frac 1{p_\oz}-1/2)]}
\frac{1}{\ro(|2^kB|)|2^kB|^{1/2}}
\|f_{k}\|_{T_2^2(\rnz)},
\end{eqnarray*}
where $2M_1>n(1/p_\oz-1/2)$. Since $f\in \twv\subset\twz$,
we have
$$\frac{1}{\ro(|2^kB|)|2^kB|^{1/2}}\|f_{k}\|_{T_2^2(\rnz)}\ls \|f\|_{\twz},$$
and for all fixed $k\in\cn$,
\begin{eqnarray*}\lim_{c\to 0}\sup_{B:\,r_B\le c}
\frac{\|f_{k}\|_{T_2^2(\rnz)}}{\ro(|2^kB|)|2^kB|^{1/2}}&&=
\lim_{c\to \fz}\sup_{B:\,r_B\ge c}
\frac{\|f_{k}\|_{T_2^2(\rnz)}}{\ro(|2^kB|)|2^kB|^{1/2}}\\&&=
\lim_{c\to \fz}\sup_{B:\,B\subset (B(0,c))^\com}
\frac{\|f_{k}\|_{T_2^2(\rnz)}}{\ro(|2^kB|)|2^kB|^{1/2}}=0.
\end{eqnarray*}
Thus, by Theorem \ref{t3.1}, we have
$$\|\pi_{L,M}f\|_{\bmor}\sim \|(t^2L)^{M_1}e^{-t^2L}
\pi_{L,M}f\|_{\twz}\ls \|f\|_{\twz},$$
and by the dominated convergence theorem for series,
\begin{eqnarray*}
\eta_1((t^2L)^{M_1}e^{-t^2L}
\pi_{L,M}f)&&=\lim_{c\rz 0}\sup_{B:\,r_B\le c}
\frac{1}{\ro(|B|)|B|^{1/2}}\bigg(\int_{B}|(t^2L)^{M_1}e^{-t^2L}
\pi_{L,M}f(x)|^2
\,dx\bigg)^{1/2}\\
&&\ls \sum_{k=0}^\fz 2^{-k[2M_1-n(1/p-1/2)]}\lim_{c\rz
0}\sup_{B:\,r_B\le c}\frac{\|f_{k}\|_{T_2^2(\rnz)}}{\ro(|2^kB|)|2^kB|^{1/2}}
=0.\nonumber
\end{eqnarray*}
Similarly, we have that $\eta_2((t^2L)^{M_1}e^{-t^2L}
\pi_{L,M}f)=\eta_3((t^2L)^{M_1}e^{-t^2L}
\pi_{L,M}f)=0$, and thus,
$(t^2L)^{M_1}e^{-t^2L}\pi_{L,M}f\in\twv$, which completes the
proof of Proposition \ref{p4.1}.
\end{proof}

\begin{lem}\label{l4.4} $\vmor\cap L^2(\rn)$ is
dense in $\vmor$.
\end{lem}
\begin{proof}
  Let $f\in\vmor$ and $M>\frac n2 (\frac 1{p_\oz}-\frac 12)$.
  Then Theorem \ref{t3.1} tells us that
$h\equiv (t^2L)^Me^{-t^2L}f\in\twv$. Similarly to the proof of
Proposition \ref{p4.1}, by Lemma \ref{l3.1}, there exist
$\{h_k\}_{k\in\cn}\subset T_{2,c}^2(\rnz)\subset \twv$ such that
$\|h-h_k\|_{\twz}\to 0$, as $k\to\fz$. Thus, by Proposition \ref{p4.1} (i)
and (iv), we obtain that $\pi_{L,1}h_k\in L^2(\rn)\cap \vmor$ and
\begin{equation}\label{4.2}
\|\pi_{L,1}(h-h_k)\|_{\bmor}\ls \|h-h_k\|_{\twz}\to 0,
\end{equation}
as $k\to\fz$.

Let $\az$ be an $(\oz,2,M,\ez)_{L^\ast}$-molecule. Then by the
definition of $H_{\oz, L^\ast}(\rn)$, $e^{-t^2L^\ast}\az\in \twl$,
which, together with Lemma \ref{l3.2}, the facts that $(\twl)^\ast=\twz$ (see
\cite[Proposition 4.1]{jyz}) and $(H_{\oz,L^\ast}(\rn))^\ast=\bmor$,
further implies that
\begin{eqnarray*}
\int_{\rn}f(x)\az(x)\,dx&&=C_M\iint_{\rnz}
(t^2L)^Me^{-t^2L}f(x)t^2L^\ast e^{-t^2L^\ast}\az(x)\frac{\,dx\,dt}{t}\\
&&=\lim_{k\to\fz}C_M\iint_{\rnz}
h_k(x,t)t^2L^\ast e^{-t^2L^\ast}\az(x)\frac{\,dx\,dt}{t}\\
&&=\frac {C_M}{C_1}\lim_{k\to\fz}
\int_{\rn}(\pi_{L,1}h_k(x))\az(x)\,dx=\frac {C_M}{C_1}\la \pi_{L,1}h,\az\ra.
\end{eqnarray*}
Since the set of finite combinations of molecules is dense
in $H_{\oz, L^\ast}(\rn)$, we then obtain that $f=\frac {C_M}{C_1}\pi_{L,1}h $
in $\bmor$.

Now, for each $k\in\cn$, let $f_k\equiv \frac {C_M}{C_1}\pi_{L,1}h_k $.
Then $f_k\in \vmor\cap L^2(\rn)$; moreover, by \eqref{4.2}, we have
$\|f-f_k\|_{\bmor}\to 0,$ as $k\to\fz$, which completes the  proof of Lemma
\ref{l4.4}.
\end{proof}

In what follows, the symbol $\la\cdot,\cdot\ra$ in the following theorem
means the duality between the space $\mathrm{BMO}_{\ro,\,L}(\rn)$
and the space $B_{\oz, L^\ast}(\rn)$ in the sense of Lemma \ref{l4.3}
with $L$ and $L^\ast$ exchanged.
Let us now state the main theorem of this paper.

\begin{thm}\label{t4.2}
The dual space of \,$\vmor$, $(\vmor)^\ast$, coincides with
the space $\bp$ in the following sense:

(i) For any $g\in B_{\oz,\,L^\ast}(\rn)$, define the linear functional $\ell$
by setting, for all $f\in \mathrm{VMO}_{\ro,\,L}(\rn)$,
\begin{equation}\label{4.3}
\ell(f)\equiv \la f, g\ra.
\end{equation}
Then there exists a positive constant
$C$ independent of $g$ such that $\|\ell\|_{(\mathrm{VMO}_{\ro,\,L}(\rn))^\ast}
\le C\|g\|_{B_{\oz,\,L^\ast}(\rn)}$.

(ii) Conversely, for any $\ell\in(\mathrm{VMO}_{\ro,\,L}(\rn))^\ast$,
there exists $g\in B_{\oz,\,L^\ast}(\rn)$, such that \eqref{4.3} holds and
there exists  a positive constant
$C$ independent of $\ell$ such that  $\|g\|_{B_{\oz,\,L^\ast}(\rn)}\le
C\|\ell\|_{(\mathrm{VMO}_{\ro,\,L}(\rn))^\ast}$.
\end{thm}

\begin{proof} By Lemma \ref{l4.3}, we have $(\bp)^\ast=\bmor$.
By Definition \ref{d3.2},  we see that $\vmor\subset \bmor$,
which implies that $\bp\subset (\vmor)^\ast$.

Conversely, let $M>\frac n2 (\frac 1{p_\oz}-\frac 12)$
and $\ell\in (\vmor)^\ast$.
By Proposition \ref{p4.1}, $\pi_{L,1}$ is bounded
from $\twv$ to $\vmor$, which implies that $\ell\circ\pi_{L,1}$
is a bounded linear functional on $\twv$. Thus, by
Theorem \ref{t4.1}, there exists $g\in \wtw$ such that
for all $f\in\twv$,  $\ell\circ \pi_{L,1}(f)=\la f,g\ra.$

Now, suppose that $f\in \vmor\cap L^2(\rn)$. By Theorem \ref{t3.1},
we have that $(t^2L)^{M}e^{-t^2L}f\in \twv$. Moreover, by
the proof of Lemma \ref{l4.4}, we see that
$f=\frac {C_M}{C_1}\pi_{L,1}((t^2L)^{M}e^{-t^2L}f)$ in $\bmor$. Thus,
\begin{eqnarray}\label{4.4}
\quad\ell(f)&&=\frac {C_M}{C_1}\ell\circ \pi_{L,1}((t^2L)^Me^{-t^2L}f)=
\frac {C_M}{C_1}\iint_{\rnz}(t^2L)^Me^{-t^2L}f(x)g(x,t)\frac{\,dx\,dt}{t}.
\end{eqnarray}

By \cite[Lemma 4.7]{jy1}, $T_{2,c}^2(\rnz)$ is dense in $\wtw$.
Since $g\in\wtw$, we may choose
$\{g_k\}_{k\in\cn}\subset T_{2,c}^2(\rnz)$ such that $g_k\to g$ in
$\wtw$. By Proposition \ref{p4.1} (iii), we have that
$\pi_{L^\ast,M}(g),\,\pi_{L^\ast,M}(g_k)\in \bp$ and
$$\|\pi_{L^\ast,M}(g-g_k)\|_{\bp}\ls
\|g-g_k\|_{\wtw}\to 0,$$
as $k\to \fz$.
This, together with \eqref{4.4}, Theorem \ref{t4.1},
the dominated convergence theorem and Lemma \ref{l4.3}, implies that
\begin{eqnarray}\label{4.5}
\ell(f)&&=\frac {C_M}{C_1}\lim_{k\to\fz}\iint_{\rnz}(t^2L)^Me^{-t^2L}f(x)g_k(x,t)
\frac{\,dx\,dt}{t}\\
&&=\frac {C_M}{C_1}\lim_{k\to\fz}\int_{\rn}f(x)\int_0^\fz
(t^2L^\ast)^Me^{-t^2L^\ast}(g_k(\cdot,t))(x)\frac{\,dt}{t}\,dx\nonumber\\
&&=\frac 1{C_1}\lim_{k\to\fz}\la f, \pi_{L^\ast,M}(g_k)\ra
=\frac 1{C_1}\la f, \pi_{L^\ast,M}(g)\ra.\nonumber
\end{eqnarray}
Since $\vmor\cap L^2(\rn)$ is dense in $\vmor$ (see Lemma \ref{l4.4}),
we finally obtain that \eqref{4.5} holds for all $f\in\vmor$,
and $\|\ell\|_{(\vmor)^\ast}=\frac 1{C_1}
\|\pi_{L^\ast,M}g\|_{\bp}$. In this sense,
we have that $(\vmor)^\ast\subset \bp$,
which completes the proof of Theorem \ref{t4.2}.
\end{proof}

\begin{rem}\rm
   If $\oz(t)\equiv t$ for all $t\in (0,\fz)$, then by
   Definitions \ref{d4.2} and \ref{d4.3}, we have that
  $\bp=H_{\oz,L^\ast}(\rn)=H_{L^\ast}^1(\rn)$, where
  $H_{L^\ast}^1(\rn)$ was the Hardy space introduced by Hofmann and Mayboroda
  in \cite{hm1}. By Theorem \ref{t4.2},
  we obtain that $(\vmor)^\ast=(\mathrm{VMO}_L(\rn))^\ast=H_{L^\ast}^1(\rn)$.
\end{rem}

In what follows, let $\oz(t)\equiv t$ for all $t\in (0,\fz)$.
We denote $(\oz,2,M,\ez)_{L}$-molecule simply by $(1,2,M,\ez)_{L}$-molecule.

We now compare the space $\vmo$ with the
classical space $\cmo$ introduced by Coifman
and Weiss \cite{cw}.

Recall that the space $\mathrm{BMO}(\rn)$ is defined to be
the set of all $f\in L^1_{\loc}(\rn)$ that satisfy
$$\|f\|_{\mathrm{BMO}(\rn)}\equiv \sup_{B\subset \rn}\frac {1}{|B|}
\int_B \lf|f(x)-f_B\r|\,dx<\fz,$$
where $f_B\equiv \frac{1}{|B|}\int_B f(x)\,dx$.
Then the space $\cmo$ is defined to be
the closure of $C_{\mathrm c}(\rn)$ (the set of all continuous functions
with compact support) in the norm $\|\cdot\|_{\mathrm{BMO}(\rn)}$.

\begin{prop}\label{p4.2} For all $n\in{\mathbb N}$,
the space $\cmo$ is continuously embedded in the space $\vmo$.
Moreover, when $n\ge 3$,
then there exists $L$ as in \eqref{1.2}
such that $\cmo$ is a proper subset of $\vmo$.
\end{prop}

\begin{proof} Let $M>\frac n4$ and
$\ez\in(0,\fz)$. Notice that each $(1,2,M,\ez)_{L^\ast}$-molecule
is a classical $H^1(\rn)$-molecule up to a harmless constant;
see \cite[p.\,41]{hm1} or \cite[Remark 7.1]{jy2}.
Thus, we have that $H_{L^\ast}^1(\rn)\subset H^1(\rn)$ and  for all
$f\in H_{L^\ast}^1(\rn)$, $\|f\|_{H^1(\rn)}\ls \|f\|_{H_{L^\ast}^1(\rn)}$.

Since $(H_{L^\ast}^1(\rn))^\ast=\bmo$ and
$(H^1(\rn))^\ast=\mathrm{BMO}(\rn)$, we further obtain that
$\mathrm{BMO}(\rn)\subset \bmo$, and for all $g\in \mathrm{BMO}(\rn)$,
$\|g\|_{\bmo}\ls\|g\|_{\mathrm{BMO}(\rn)}.$

Let us now show that $C_{\mathrm c}(\rn)\subset \vmo$.
Suppose that $f\in C_{\mathrm c}(\rn)$  and $B\equiv B(x_B,r_B)$.
By an argument as in \cite[Section 2.5]{a1},
we have that for all $t>0$, $e^{-tL}1=1$ in the
$L_{\loc}^2(\rn)$ sense, that is, for any $\phi\in L^2(\rn)$
with compact support, there holds
$$\int_{\rn}\phi(x)\,dx=\int_{\rn}e^{-tL}1 \phi(x)\,dx
=\int_{\rn}e^{-tL^\ast}\phi(x)\,dx.$$
This together with Lemma \ref{l2.2} and the H\"older inequality implies that
\begin{eqnarray*}
 \|(I-e^{-r_B^2L})^Mf\|_{L^2(B)}&&=
 \sup_{\|g\|_{L^2(B)}\le 1}\lf|\int_{B}(I-e^{-r_B^2L})^Mf(x) g(x)\,dx\r|\\
 &&=\sup_{\|g\|_{L^2(B)}\le 1}\lf|\int_{B}(I-e^{-r_B^2L})^M(f-f_B)(x) g(x)\,dx\r|\\
 &&\ls \sup_{\|g\|_{L^2(B)}\le 1}\sum_{k=0}^\fz \exp\lf\{-\frac{2^kr_B^2}{cr_B^2}\r\}
 \|f-f_B\|_{L^2(U_k(B))}\|g\|_{L^2(B)}\\
 &&\ls \sum_{k=0}^\fz e^{-2^k/{c}} \lf\{\|f-f_{2^kB}\|_{L^2(U_k(B))}
 +|U_k(B)|^{1/2}|f_B-f_{2^kB}|\r\}\\
 &&\ls \sum_{k=0}^\fz e^{-2^k/{c}}2^{kn/2}\|f-f_{2^kB}\|_{L^2(U_k(B))},
\end{eqnarray*}
where $c$ is a positive constant.
Since $f\in C_{\mathrm c}(\rn)\subset \mathrm{CMO}(\rn)$,
by the dominated convergence theorem for series,
\cite[Proposition 3.5]{ddsty} and the fact
that $\mathrm{CMO}(\rn) \subset
\mathrm{VMO}(\rn)$, we have
\begin{eqnarray*}
\gz_1(f)&&=\lim_{c\to0}\sup_{B:\,r_B\le c}
\frac{1}{|B|^{1/2}}\|(I-e^{-r_B^2L})^Mf\|_{L^2(B)}\\
&&\ls
\sum_{k=0}^\fz e^{-2^k/{c}}2^{kn}\lim_{c\to0}
\sup_{B:\,r_B\le c}\frac{1}{|2^kB|^{1/2}}
\|f-f_{2^kB}\|_{L^2(U_k(B))}=0.
\end{eqnarray*}
Similarly, we have that $\gz_2(f)=\gz_3(f)=0$, which implies that
$f\in\vmo$, and hence $C_{\mathrm c}(\rn)\subset \vmo$.
Since $C_{\mathrm c}(\rn)$ is dense in $\cmo$, by the fact that
for all $g\in \mathrm{CMO}(\rn)$,
$\|g\|_{\vmo}=\|g\|_{\bmo}\ls\|g\|_{\mathrm{BMO}(\rn)}
\sim\|g\|_{\mathrm{CMO}(\rn)},$ we finally obtain that $\cmo\subset \vmo$.

On the other hand, if $n\ge 3$, then there exist $L$
as in \eqref{1.2} and $p\in (1,2)$ such that
$\nabla (L^{\ast})^{-1/2}$ is not bounded on $L^p(\rn)$; see \cite{a1,f}.
By the fact that the Riesz transform $\nabla (L^{\ast})^{-1/2}$ is
bounded on $L^2(\rn)$ and from $H_{L^\ast}^1(\rn)$ to $L^1(\rn)$
(see \cite{hm1,jy2}), we obtain that $H_{L^\ast}^1(\rn)$
is a proper subspace of $H^1(\rn)$. Otherwise, $H_{L^\ast}^1(\rn)=H^1(\rn)$,
then by the interpolation theorem, $\nabla (L^{\ast})^{-1/2}$ is
bounded on $L^p(\rn)$ for all $p\in (1,2)$,
which contradicts with the fact that
$\nabla (L^{\ast})^{-1/2}$ is not bounded on $L^p(\rn)$ for some
$p\in (1,2)$.

Now by Theorem \ref{t4.2}, we have that
$$({\mathrm{VMO}_L(\rn)})^\ast=H_{L^\ast}^1(\rn)\subsetneqq H^1(\rn)=(\cmo)^\ast,$$
which implies that $\cmo\subsetneqq \vmo$, and hence, completes the
proof of Proposition \ref{p4.2}.
\end{proof}

From Proposition \ref{p4.2}, we deduce the following conclusion.

\begin{cor}\label{c4.1}
When $A$ has real entries, or when the dimension $n=1$ or $n=2$ in
the case of complex entries, the spaces $\cmo$ and $\vmo$ coincide with
equivalent norms.
\end{cor}

To prove Corollary \ref{c4.1}, we need the following
lemma. (We are very grateful to the referee for him/her
to tell us Corollary \ref{c4.1} and its proof.)

\begin{lem}\label{l4.5}
Let $X$ and $Y$ be two Banach spaces such that $X\subset Y$
with continuous embedding. Assume that $X^\ast$ and $Y^\ast$
coincide with equivalent norms.
Then $X= Y$ with equivalent norms.
\end{lem}

\begin{proof} Let $J : X\to Y$ be the inclusion map, which is
a continuous linear map by assumption. As is easily seen,
its adjoint $J^\ast:\ Y^\ast\to X^\ast$ is the map which, to any $\phi\in  Y^\ast$,
associates its restriction to $X$. The Hahn-Banach theorem yields
that $J^\ast$ is onto. Moreover, it is also one-to-one. Indeed,
if this was not true, there would exist $\phi\in Y^\ast$ identically vanishing
on $X$ but not on $Y$, which would contradict the assumption
that $X^\ast$ and $Y^\ast$ coincide with equivalent norms. Thus,
$J^\ast$ is a continuous isomorphism, and, by Theorem 4.15 in \cite{ru},
$J$ is an isomorphism between $X$ and $Y$, which exactly means that
$X=Y$. From this and Corollaries 2.12(c) of the open mapping
theorem in \cite[pp.\,49-50]{ru}, it follows that $X$
and $Y$ coincide with equivalent norms, which completes the proof
of Lemma \ref{l4.5}.
\end{proof}

\begin{proof}[Proof of Corollary \ref{c4.1}]
When $A$ has real entries, or when the dimension $n=1$ or $n=2$ in
the case of complex entries, the heat kernel always satisfies the
Gaussian pointwise estimates in size and regularity (see \cite{at}),
and the spaces $H_{L^\ast}^1(\rn)$ and $H^1(\rn)$ coincide
with equivalent norms (see \cite{ar,ya08}).
By this and Proposition \ref{p4.2}, we obtain that $\cmo\subset\vmo$
and $(\cmo)^\ast=H^1(\rn)=H_{L^\ast}^1(\rn)=(\vmo)^\ast$,
which together with Lemma \ref{l4.5} imply that $\cmo$ and $\vmo$ coincide with
equivalent norms. This finishes the proof of Corollary \ref{c4.1}.
\end{proof}

{\bf\noindent Acknowledgment}

\medskip

The authors sincerely wish to express their deeply thanks to the
referee for her/his very carefully reading and also her/his so
many valuable and helpful remarks which made this article more
readable. In particular, Corollary \ref{c4.1} and its proof belong
to the referee.

\medskip

\noindent Renjin Jiang:

\noindent School of Mathematical Sciences, Beijing Normal University

\noindent Laboratory of Mathematics and Complex Systems, Ministry
of Education

\noindent 100875 Beijing

\noindent People's Republic of China

\noindent E-mail: \texttt{rj-jiang@mail.bnu.edu.cn}

\medskip

\noindent {\it Present address}

\noindent Department of Mathematics and Statistics,
University of Jyv\"{a}skyl\"{a}

\noindent P.O. Box 35 (MaD), 40014, Finland

\medskip
\noindent Dachun Yang (corresponding author)

\noindent School of Mathematical Sciences, Beijing Normal University

\noindent Laboratory of Mathematics and Complex Systems, Ministry
of Education

\noindent 100875 Beijing

\noindent People's Republic of China

\noindent E-mail: \texttt{dcyang@bnu.edu.cn}

\end{document}